\documentclass[12pt,centertags,oneside]{amsart}

\pdfoutput=1

\usepackage{tikz}
\usetikzlibrary{shapes.geometric}
\usetikzlibrary {angles,calc,quotes}

\usepackage{circuitikz}
 \usepackage{caption}
 \usepackage{amsmath}
\usepackage{amsmath,amstext,amsthm,amscd,typearea,hyperref,stmaryrd,mathabx}
\usepackage{amssymb}
\usepackage{a4wide}
\usepackage[mathscr]{eucal}
\usepackage{mathrsfs}
\usepackage{typearea}
\usepackage{charter}
\usepackage{pdfsync}
\usepackage[a4paper,width=16.2cm,top=3cm,bottom=3cm,lines=50]{geometry}

\numberwithin{equation}{section}

\newtheorem{theorem}{Theorem}[section]
\newtheorem{definition}[theorem]{Definition}
\newtheorem{proposition}[theorem]{Proposition}
\newtheorem{corollary}[theorem]{Corollary}

\newtheorem{lemma}[theorem]{Lemma}
\newtheorem{remark}[theorem]{Remark}

\newcommand{\cali}[1]{\mathscr{#1}}

\newcommand{\PGL}{{\rm PGL}}

\newcommand{\diag}{{\rm diag}}

\newcommand{\Tan}{{\rm Tan}}

\newcommand{\sing}{{\rm sing}}

\newcommand{\rank}{{\rm rank}}

\renewcommand{\Im}{{\rm Im}}
\newcommand{\nd}{{\rm nd}}

\newcommand{\Ac}{\cali{A}}

\newcommand{\Cc}{\cali{C}}

\newcommand{\Fc}{\cali{F}}
\newcommand{\Gc}{\cali{G}}

\newcommand{\Oc}{\cali{O}}

\newcommand{\Uc}{\cali{U}}

\newcommand{\Sc}{\cali{S}}

\newcommand{\C}{\mathbb{C}}

\renewcommand{\H}{\mathbb{H}}
\newcommand{\N}{\mathbb{N}}

\newcommand{\R}{\mathbb{R}}

\newcommand{\U}{\mathbb{U}}

\renewcommand{\P}{\mathbb{P}}

\newcommand{\G}{\mathbb{G}}

\newcommand{\al}{\alpha}

\newcommand{\rvline}{\hspace*{-\arraycolsep}\vline\hspace*{-\arraycolsep}}

\title[]{Generic singularities of  holomorphic foliations by curves on $\P^n$}

\author{Sahil Gehlawat}
\address{Universit\'e de Lille, 
Laboratoire de math\'ematiques Paul Painlev\'e, 
CNRS U.M.R. 8524,  
59655 Villeneuve d'Ascq Cedex, 
France.} 
\email{sahil.gehlawat@univ-lille.fr}

\author{Vi{\^e}t-Anh Nguy{\^e}n}
\address{Universit\'e de Lille, 
Laboratoire de math\'ematiques Paul Painlev\'e, 
CNRS U.M.R. 8524,  
59655 Villeneuve d'Ascq Cedex, 
France.} 
\email{Viet-Anh.Nguyen@univ-lille.fr, {\tt  https://pro.univ-lille.fr/viet-anh-nguyen/}}

\address{and Vietnam Institute for Advanced Study in Mathematics (VIASM),  157 Chua Lang Street, Hanoi, Vietnam.
}	

\date{September 09, 2024}

\begin{document}

%for line number
%\internallinenumbers

\begin{abstract}
Let   $\Fc_d(\P^n)$ be the space of all singular holomorphic  foliations  by curves on $\P^n$ ($n\geq 2$) with degree $d\geq 1.$
We  show  that there is     subset  $\Sc_d(\P^n)$ of  $\Fc_d(\P^n)$  with full  Lebesgue  measure with the following  properties:
\begin{itemize}  \item[$\bullet$] for every $\Fc\in\Sc_d(\P^n),$
 all singular points of  $\Fc$ are linearizable hyperbolic.
 \item[$\bullet$]
 If, moreover, $d\geq 2,$ then every   $\Fc$ does not possess any invariant algebraic  curve.
 \end{itemize}
\end{abstract}

%\medskip\medskip

%\noindent
%\footnote{{\bf MSC 2020:} Primary: 37F75,  37A30;  Secondary: 57R30.}
\subjclass[2020]{Primary: 37F75,  37A30;  Secondary: 57R30}
% \medskip

%\noindent
%{\bf Keywords:} singular  holomorphic  foliation, hyperbolic singularity, linearizable singularity, invariant algebraic  %curve.
\keywords{singular  holomorphic  foliation, hyperbolic singularity, linearizable singularity, invariant algebraic  curve}

\maketitle

%%%%%%%%%%%%%%%%%%%%%%%%%%%%%%%%%%%%%%%%%%%%%%%%%%%%%%%%%%%%%%%%%%%%%%%%%%%%%%
%%%%%%%%%%%%%%%%%%%%%%%%%%%%%%%%%%%%%%%%%%%%%%%%%%%%%%%%%%%%%%%%%%%%%%%%%%%%%%
\section{Introduction} \label{S:Intro}
%%%%%%%%%%%%%%%%%%%%%%%%%%%%%%%%%%%%%%%%%%%%%%%%%%%%%%%%%%%%%%%%%%%%%%%%%%%%%%
%%%%%%%%%%%%%%%%%%%%%%%%%%%%%%%%%%%%%%%%%%%%%%%%%%%%%%%%%%%%%%%%%%%%%%%%%%%%%%

 Let $n\in\N$ with $n\geq 2.$ A  (singular) holomorphic foliation  (by curves) on $\P^n$  with degree $d\in\N$  is given by  a bundle  morphism
$$
\Phi:\ \Oc_{\P^n}(1-d)\to\Tan(\P^n),
$$
 and 
its singular set is 
$$
\sing(\Fc):=\{x\in\P^n:\ \Phi(x)=0\}.
$$
Here  $\Oc_{\P^n}(-d):= (\Oc_{\P^n}(-1))^{\otimes d},$ where $\Oc_{\P^n}(-1)$ is the tautological line bundle on $\P^n, $ and     $\Tan(\P^n)$ is the  tangent bundle of $\P^n.$
  A   holomorphic  foliation $\Fc$ of degree $d\in\N$ on $\P^n$ is  represented  in affine  coordinates $x=(x_1,\ldots,x_n)$  by a vector field  of the  form 
  \begin{equation}\label{e:representation}
   X=\sum_{j=0}^{d} R_j+ g\mathcal R,
  \end{equation}
where $\mathcal R$ is the radial  vector field  defined by $\mathcal R(x):= \sum_{j=1}^n x_j{\partial\over \partial x_j},$ $g$ is a homogeneous  
polynomial of degree $d,$ and $R_j$ is a vector   field whose components  are  homogeneous polynomial of degree $j,$
$0\leq j\leq d.$
Since $\sing(\Fc)$ has
codimension greater than 1 we have either $g\not\equiv  0$ or $ g \equiv 0$ and $R_d$ cannot
be written as $hR$ where $h$ is homogeneous of degree $d -1.$ In this case
$X$ has a pole of order $d - 1$ at infinity (see, for example, \cite{LinsNetoScardua}). We call $d$ the degree of
the foliation.
Let $\Fc_d(\P^n)$ denote  the (moduli) space of  singular  holomorphic  foliations by curves   of degree $d$ on $\P^n.$
Using representation \eqref{e:representation}, $\Fc_d(\P^n)$ is  identified to a Zariski open set   of the projective  space $\P^M$  with
\begin{equation}\label{e:M}
M=M(n,d):= n {n+d\choose d}  + {n+d-1\choose  d}-1. 
\end{equation}
The space $\Fc_d(\P^n)$  plays an important role in the theory of singular holomorphic  foliations by curves  as it  is  one of the main  sources  of examples,  and  hence it is  often  used  in the research  as  a testing ground for expecting results.
Moreover, 
since every holomorphic  foliation $\Fc$ on $\P^n$ is  singular, i.e.,  $\sing(\Fc)\not=\varnothing,$  it  is   natural  to investigate  the  singular  set  $\sing(\Fc).$  At the first stage,  one may    understand
the singular set   $\sing(\Fc)$ of a generic element $\Fc\in\Fc_d(\P^n),$ that is, $\Fc\in\Uc,$ where $\Uc\subset \Fc_d(\P^n)$ such that the complement  $\Fc_d(\P^n)\setminus \Uc$ is,  in some  sense, negligible. Another basic question is  to investigate if a 
generic element  $\Fc\in\Fc_d(\P^n)$  admits an invariant algebraic curve.

The pioneering  article of Jouanolou \cite{Jouanolou} is the first work to  study this problem systematically.
Jouanolou  obtains a satisfactory description  of a generic foliation with a given degree $d\geq 1$  on $\P^2.$  One of the main ingredients
in his  method is  to introduce  an explicit family of foliations $\mathcal J^{2,d}$ on $\P^2$  which enjoys 
some remarkable properties.  Another ingredient in his proof is
a perturbation  argument which  roughly says that if one foliation  satisfies some  properties, then every
nearby foliation (in the moduli space of foliations)
should do so as well. Using  the holomorphic propagation, this argument  permits to show that every  foliation
$\Fc$ in a subset $\Uc$  of full Lebesgue measure of $\Fc_d(\P^n)$ satisfies these properties. Later on Lins Neto 
\cite{LinsNeto} improves the result of Jouanolou by showing that $\Uc$ can be chosen to be  open.
Soares  in \cite{Soares} proves   Jouanolou's result for $\P^3.$

The family of Jouanolou foliations was later  generalized to $\mathcal J^{n,d}$ by  Lins Neto and Soares in \cite{LinsNetoSoares,Soares}
for  general projective spaces $\P^n,$ $n\geq 2.$ These  two  authors develop a residue calculus which allows them to detect effectively if a given foliation  has an invariant  algebraic curve.

The  following  theorem    summarizes  the  above mentioned results     describing  some typical properties of a generic foliation $\Fc\in\Fc_d(\P^n).$ 

Set 
\begin{equation}\label{e:N}
N=N(n,d):=1+d+\ldots+d^n.
\end{equation}
 \begin{theorem}\label{T:generic-JNS} 
  {\rm (Jouanolou \cite{Jouanolou} and Lins Neto \cite{LinsNeto} for $n=2,$  Soares \cite{Soares} for $n=3,$ Lins Neto--Soares \cite{LinsNetoSoares} for $n\geq 2$).}  For every $n,d\in\N$ with $n\geq 2$ and $d\geq 2,$
  there is a nonempty real Zariski  open set $\Uc^{\nd}_d(\P^n)\subset  \Fc_d(\P^n)$ such that for every $\Fc\in \Uc^{\nd}_d(\P^n),$
  all  the singularities of $\Fc$ are    non-degenerate and  $\Fc$ does not possess any invariant algebraic  curve. Moreover, every $\Fc\in \Uc^{nd}_d(\P^n)$  admits exactly  $N(n,d) $  singular points.
\end{theorem}

It is  worth noting the   work of  Loray--Rebelo \cite{LorayRebelo} investigating, among other questions,  whether    a generic foliation $\Fc\in \Fc_d(\P^n)$ admits an invariant  algebraic  subvariety of dimension $\geq 1$ (see Remark  \ref{R:LorayRebelo}  below).
  
It turns out that  foliations  with only non-degenerate singularities may not be the right   context  for many problems  arising from the theory of singular holomorphic  foliations.
Indeed,  one of the main developments in this field during  the last 25 years is the advent of the  ergodic  theory which emphasizes a global  study of singular holomorphic foliations
using     the analysis of positive currents and the  tools from  geometric complex    analysis developed by     Berndtsson--Sibony \cite{BerndtssonSibony}  and Forn{\ae}ss--Sibony \cite{FornaessSibony05,FornaessSibony08,FornaessSibony10}.
One major goal of the latter theory is to  
study the  global behaviour of  a  generic  leaf of a foliation $\Fc.$  Here, the  term {\it ``generic''} is measured  with respect to   some new objects called {\it  directed positive  harmonic  currents} (see for example \cite[Definition 7]{FornaessSibony08} etc. for a definition of this  notion).  
The      analysis  of these  currents 
are often  difficult when one approaches  the singularities of  $\Fc.$
The study near  a  singular point will be more effective, more tangible  if  it is   linearizable hyperbolic  singular point.        
Indeed, in this case it admits (local) invariant coordinate hyperplanes, and  hence it   allows us  to start the  analysis of the   directed positive  harmonic  currents of $\Fc$ from   these special  coordinates, and consequently  results in quite  a lot of important geometric, dynamical properties such as the regularity of  the leafwise Poincar\'e metric on $\Fc,$ the  entropies, the Lyapunov exponents,     the unique ergodicity etc. in the works of Dinh, Forn{\ae}ss, Sibony and  many others  (see \cite{DinhNguyenSibony12, DinhNguyenSibony14,  DinhNguyenSibony22, DinhSibony18,DinhWu,  FornaessSibony10,NguyenVietAnh18a,NguyenVietAnh18b, NguyenVietAnh23,  NguyenVietAnh24} etc. and the references therein, see also the  surveys \cite{DinhSibony20,FornaessSibony08,NguyenVietAnh18c,NguyenVietAnh20b}).  
When  the singular points of $\Fc$  are not all   hyperbolic,  some  of these desirable properties may not hold as illustrated  by Chen \cite{Chen} and 
Alkateeb--Rebelo  \cite{AlkateebRebelo}. 
When  the  singular points of $\Fc$ are  merely      
non-degenerate (and not linearizable), the analysis becomes much  harder, and only  a few    results have been  obtained,  see the articles of Bacher \cite{Bacher23, Bacher24} etc.  which  are  in fact  inspired  by the corresponding results in   the  case of  linearizable singular points.
   
   The main purpose of  this article is  to prove   Theorem \ref{T:generic-JNS} in the  context
   of  linearizable hyperbolic  singular points. Here is  our main result.
  
\begin{theorem}\label{T:main} {\rm (Main Theorem)} For every $n,d\in\N$ with $n\geq 2$ and $d\geq 1,$
there is a subset $\Sc_d(\P^n)$ of full Lebesgue measure of $\Fc_d(\P^n)$ such that for every $\Fc\in\Sc_d(\P^n),$ all singular points of $\Fc$ are   linearizable hyperbolic  and  that  if $d\geq 2,$ then $\Fc$ does not possess any invariant algebraic  curve.  
 \end{theorem}
  
  The case $n=2$ of Theorem \ref{T:main}  has   already been known in the   work of Jouanolou  \cite{Jouanolou}.  Theorem \ref{T:main} provides an abundant  supply of foliations   for which  the ergodic theory of  singular holomorphic foliations   can be explored efficiently by means of  the  theory of positive  currents and the theory of  geometric  complex  analysis.

  Theorem  \ref{T:main}, combined with  Corollary 1.3  in \cite{NguyenVietAnh24}, implies the following result:
\begin{corollary}\label{C:Nguyen}
 Let $\Fc \in \Sc_d(\P^n).$
Then for every  positive  harmonic current $T$  directed by $\Fc,$  $T$ is  diffuse and  the   Lelong  number of $T$ vanishes everywhere   in $\P^n.$ In particular,  this  property holds for every foliation in a subset of full Lebesgue measure of $\Fc_d(\P^n).$
\end{corollary}   
  
  At the first glance,  the conclusion of Theorem \ref{T:main} appears to be intuitively evident  since in the local  situation,  for a generic (in the sense of Lebesgue measure)   polynomial vector field $Z$ singular at $0\in\C^n,$  the point $\{0\}$ is linearizable  hyperbolic. 
  However,  the  genericity in the moduli space $\Fc_d(\P^n)$  is  not  such evident at all as  it is now  a global  situation.
  The  main  ingredient in the proof of   Theorem \ref{T:main} is then  the  introduction of a  new $n$-dimensional family $\{\mathcal J_\alpha,$ $\alpha\in\C^n\}$
  of singular holomorphic foliations having the Jouanolou foliation  $\mathcal J^{n,d}$ as $\mathcal J_0$ (see Definition \ref{D:new-family} below). This family has  the  novelty  that the  characteristic polynomial  of every  singular point of the foliation $\mathcal J_\alpha$ varies biholomorphically as the parameter  $\alpha$ varies in  a small open neighborhood $\U$    of $0\in\C^n$ (see Theorem \ref{T:submersion}). This  novelty will make  the previous local intuition  come true in the global context.  We refer the readers to 
  Remark  \ref{R:LinsNetoSoares} below for a discussion of the choice of another family
  close to that   introduced by Lins Neto--Soares \cite{LinsNetoSoares}.
  
  We   use  the method   developed  by   Lins Neto--Soares \cite{LinsNetoSoares} in order to show that
  for almost every $\alpha\in\U, $ the foliation $\mathcal J_\alpha$ does not possess an invariant algebraic  curve.
  Finally,  we conclude the proof using  an argument of holomorphic propagation.
  
  \smallskip
  
  The  article is organized as  follows.
  
  \smallskip
  
  In Section \ref{S:Background}
we recall  basic notions and results on  the local classification of singularities of holomorphic vector fields. We also  recall some basic properties of the Jouanolou foliations.   
  We start Section \ref{S:New-family}  by introducing  the new family of foliations  $\{\mathcal J_\alpha:\ \alpha\in\C^n\}.$ After a qualitative  study of  their singularities,  the section deals   with the proof of its main result Theorem \ref{T:submersion}  which is  one of the main tools of our  method.
Section   \ref{S:Main-Theorem } is devoted   to the proof of the Main Theorem \ref{T:main}.
 More concretely,  Subsection  \ref{SS:No-invariant-curve}  shows that  $\Fc_\alpha$
has no invariant algebraic curve for a generic value of $\alpha\in\U.$  Combining the results of  Section \ref{S:New-family} and  Subsection  \ref{SS:No-invariant-curve}, Subsection \ref{SS:End}            concludes  the  proof.  The article is ended with some remarks.

 \smallskip
 
\noindent
{\bf Acknowledgments. }  The  authors   acknowledge support from the Labex CEMPI (ANR-11-LABX-0007-01)
and from  the project QuaSiDy (ANR-21-CE40-0016).
The paper was partially prepared 
during the visit of the second  author at the Vietnam  Institute for Advanced Study in Mathematics (VIASM). He would like to express his gratitude to this organization 
for hospitality and  for  financial support.

%%%%%%%%%%%%%%%%%%%%%%%%%%%%%%%%%%%%%%%%%%%%%%%%%%%%%%%%%%%%%%%%%%%%%%%%
%%%%%%%%%%%%%%%%%%%%%%%%%%%%%%%%%%%%%%%%%%%%%%%%%%%%%%%%%%%%%%%%%%%%%%%%
\section{Background}\label{S:Background}
%%%%%%%%%%%%%%%%%%%%%%%%%%%%%%%%%%%%%%%%%%%%%%%%%%%%%%%%%%%%%%%%%%%%%%%%%
%%%%%%%%%%%%%%%%%%%%%%%%%%%%%%%%%%%%%%%%%%%%%%%%%%%%%%%%%%%%%%%%%%%%%%%%%

 We first recall some basic elements of    holomorphic  vector fields and  the local theory of singular holomorphic foliations.
 \begin{definition}
\label{D:singularities}\rm
  Let $Z=\sum_{j=1}^n f_j(z){\partial\over\partial  z_j}$ be  a holomorphic vector field  defined in a neighborhood $\U$ of $0\in\C^n.$ Consider  the holomorphic map
  $f:=(f_1,\ldots,f_n):\  \U\to\C^n.$ We say that $Z$ is 
\begin{enumerate}
\item {\it  singular  at $0$} if  $f(0)=0.$
%\item {\it  linear} if
%it can be written as    $$\sum_{j=1}^n(\sum_{\ell=1}^n a_{j\ell}z_\ell)  {\partial \over \partial z_j}\quad\text{for some %constants}\quad  a_{j\ell}\in\C,\quad 1\leq j,\ell\leq  n.$$
\item {\it generic linear} if
it can be written as 
$$Z(z)=\sum_{j=1}^n \lambda_j z_j {\partial\over \partial z_j}$$
where $\lambda_j$ are non-zero complex numbers. The  $n$  hyperplanes $\{z_j=0\}$ for $1\leq j\leq n$ are said to be the {\it invariant  hypersurfaces.}

\item     {\it with non-degenerate singularity at $0$} if $Z$ is  singular at $0$ and the eigenvalues $\lambda_1,\ldots,\lambda_n$ of the Jacobian matrix $Df(0)$
 are all  nonzero. 
 %We say  that the singularity is in the Poincar\'e domain if the convex hull in $\C$ of
%$\{\lambda_1 , \ldots , \lambda_n \}$ does not contain the origin, it is in the Siegel domain otherwise.

%\item    {\it with weakly hyperbolic  singularity at $0$} if $Z$ is  singular at $0$ and the eigenvalues $\lambda_1,\ldots,%\lambda_n$ of the Jacobian matrix $Df(0)$
 %are all  nonzero and  there  are  some $1\leq j\not=l\leq k$ with $\lambda_j/\lambda_l\not\in\R.$   

\item    {\it with hyperbolic  singularity at $0$} if $Z$ is  singular at $0$ and the eigenvalues $\lambda_1,\ldots,\lambda_n$ of the Jacobian matrix $Df(0)$
 are all  nonzero and $\lambda_j/\lambda_l\not\in\R$  for all $1\leq j\not=l\leq n.$
\end{enumerate}
The integral curves of $Z$ define a  singular holomorphic foliation on
$\U$.
The condition $\lambda_j\not=0$ implies that the foliation   
has an isolated singularity at 0.
The singular set of $Z$ is denoted by $\sing(Z).$
\end{definition}

\begin{definition}\rm 
 Let $\Fc$ be a singular holomorphic foliation  such that  its  singular set $\sing(\Fc)$ is an  analytic subset of $X$  with codimension $\geq 2.$  
 We say that a
singular point $a\in \sing(\Fc) $ is {\it linearizable} (resp. {\it non-degenerate},   resp.  {\it hyperbolic}) if 
there is a local holomorphic coordinate system of $X$ near $a$ on which the
leaves of $\Fc$ are integral curves of a generic linear vector field  (resp.  of a holomorphic vector field admitting $0$ as non-degenerate singularity,   resp. a  hyperbolic singularity).
\end{definition}

\begin{definition}
 \label{D:nu} \rm  Given  $\delta>0,$ a point $\lambda=(\lambda_1,\ldots,\lambda_n)\in\C^n$  is  said to be   of type $\delta$  if there is  a constant $c>0$ such that for every $1\leq j\leq n,$ we have
 $$
 |\lambda_j-\sum_{i=1}^n m_i\lambda_i|\geq  {c\over  |m|^\delta},
 $$
 for all vectors $m=(m_1,\ldots,m_n)\in\N^n$ with $|m|:=\sum_{i=1}^n m_i\geq 2.$
\end{definition}

\begin{proposition}\label{P:Arnold}
 For every $\delta>{n-2\over n},$ consider the set
 $$
 \Ac_\delta:= \left\lbrace  \lambda\in\C^n:\  \lambda \text{  of type $\delta$ } \right\rbrace.
 $$
 Then the $2n$-Lebesgue measure of the set $\C^n\setminus \Ac_\delta$ is  zero.
\end{proposition}
\proof
It is an  immediate consequence of   \cite[Theorem 4, p. 189]{Arnold}.
\endproof

\begin{theorem}\label{T:Siegel} {\rm  (Siegel's Theorem, see for example \cite{IY})} If  the eigenvalues of the linear part of a holomorphic vector  field at a singular point form  a vector of type  $\delta,$ then the  field is  biholomorphically   equivalent  to its linear part in the neighborhood of  the  singular point. 
\end{theorem}

In the  remainder of the section we  state  some useful  results in the  global context of singular holomorphic foliations on $\P^n.$ 
\begin{definition} \label{D:F-nd}\rm 
 For  $n\geq 2$ and $d\geq 1,$ let $\Fc^{\nd}_d(\P^n)$ denote the space of all foliations $\Fc\in\Fc_d(\P^n)$ 
 such that all  singular points of $\Fc$  are  non-degenerate.
\end{definition}

\begin{proposition}\label{P:Local-holo-charts}{\rm  (See for example \cite[Corollary 2.7]{LinsNetoScardua}) }  
 Let $\Fc_0\in \Fc^{\nd}_d(\P^n)$ with  $\sing(\Fc_0)=\{p_{1,0},\ldots, p_{N,0}\},$ where $p_{i,0}\not=p_{j,0}$ for $i\not=j.$  Then  there exist  connected open neighborhood $\mathcal U$ of $\Fc_0$ in $\Fc_d(\P^n),$   $U_1,\ldots, U_N$
 of $p_{1,0},\ldots, p_{N,0}$ respectively,  pairwise disjoint, and holomorphic maps $\varphi_j:\ \mathcal U\to  U_j,$
 for $1\leq j\leq N$  satisfying  the  following properties:
 
\begin{itemize}
 \item [(i)] $\varphi_j(\Fc_0)=p_{j,0}.$
 \item[(ii)] For every   $\Fc\in\mathcal U$ and every $1\leq j\leq N,$ $\varphi_j(\Fc)$ is the unique singularity of $\Fc$ in $U_j,$ and  moreover,   this singularity is   non-degenerate.
 \item[(iii)]  For every $\Fc\in\mathcal U,$ we have  $\sing(\Fc)=\{ \varphi_1(\Fc),\ldots,\varphi_N(\Fc)  \}.$
\end{itemize}
In particular,  $ \Fc^{\nd}_d(\P^n)$ is a Zariski  open subset  of $\Fc_d(\P^n),$ and hence it is  connected.
\end{proposition}

Consider the Jouanolou foliation $\mathcal{J}^{n,d} \in \Fc_d(\P^n)$ given by the vector field (see \cite{Jouanolou,Soares,LinsNetoSoares}):
\begin{equation}\label{e:Jouanolou}
X_0 = \sum_{i=1}^{n-1}(x_{i+1}^d - x_{i}x_1^{d})\frac{\partial}{\partial x_i} + (1 - x_n x_1^{d})\frac{\partial}{\partial x_n}.
\end{equation}

\begin{proposition}\label{P:Jouanolou}{\rm  (See for example \cite[Proposition  2.10]{LinsNetoScardua}) } 
 \begin{enumerate}
  \item $\mathcal{J}^{n,d} \in \Fc^{\nd}_d(\P^n).$
  \item $\sing(\mathcal{J}^{n,d})\subset \C^n.$ Moreover,  
  $\sing(\mathcal{J}^{n,d})=\{p_1,\ldots,  p_N\},$ where $N=N(n,d)$ and 
  $
  p_j=(\xi^j, (\xi^j)^{f(n-1)}, \ldots,  (\xi^j)^{f(1)}),$  $j=1,\ldots,N,$ where 
  $\xi$ is an $N$-th primitive root of  unity and $f(m):=-(d+d^2+\ldots+d^m).$
 \end{enumerate}
\end{proposition}

%%%%%%%%%%%%%%%%%%%%%%%%%%%%%%%%%%%%%%%%%%%%%%%%%%%%%%%%%
%%%%%%%%%%%%%%%%%%%%%%%%%%%%%%%%%%%%%%%%%%%%%%%%%%%%%%%%%%
\section{A new family of foliations} \label{S:New-family}
%%%%%%%%%%%%%%%%%%%%%%%%%%%%%%%%%%%%%%%%%%%%%%%%%%%%%%%%%%
%%%%%%%%%%%%%%%%%%%%%%%%%%%%%%%%%%%%%%%%%%%%%%%%%%%%%%%%%%
\begin{definition}\label{D:new-family}\rm  For $\alpha = (\alpha_1, \alpha_2, \ldots, \alpha_n) \in \mathbb{C}^n$, consider the singular holomorphic foliation $\mathcal J_\alpha$ on $\P^n$  which is generated on $\C^n$ by the  vector field
\begin{equation}\label{e:Jouanolou-new}
X_{\alpha} = \sum_{i=1}^{n-1}(x_{i+1}^d - x_{i}x_1^{d} + \alpha_{i})\frac{\partial}{\partial x_i} + (1 - x_n x_1^{d} + \alpha_{n})\frac{\partial}{\partial x_n} = X_{0} + \sum_{i=1}^{n} \alpha_{i} \frac{\partial}{\partial x_{i}}.
\end{equation}
\end{definition}
 Given  an open  neighourhood $\U$ of $0 \in \mathbb{C}^n$, the set $\mathcal{J}(\U) := \{\mathcal{J}_{\alpha} \in \Fc_d(\P^n) : \alpha \in \U\}$ is an $n$-dimensional subspace in $\Fc_d(\P^n)$ passing through $\mathcal J_0=  {\mathcal J}^{n,d}$.  By  Proposition \ref{P:Jouanolou},  $\mathcal J_0$ has   $N=N(n,d) $ number of singular points, each of which is non-degenerate. Choose $\U$ to be small enough so that, each $\mathcal J_{\alpha}$ for $\alpha \in \U$ also has exactly $N $  singular points and that the set  $E_{\alpha}:=\sing(X_{\alpha})$ coincides with $\sing(\mathcal J_\alpha).$ 
Observe  that  $x=(x_{1,\alpha}, x_{2,\alpha}, \ldots x_{n,\al}) \in E_{\alpha}$ if and only if
\begin{equation}
 \label{e:E_alpha}
 x_{i+1}^d - x_{i}x_1^{d} + \alpha_{i}=0\quad\text{for}\quad 1\leq i\leq n-1,\quad\text{and}\quad 1 - x_n x_1^{d} + \alpha_{n}=0.
\end{equation}
For every $1 \le m \le N$, let \begin{equation}\label{e:p_m} p_{m} = (\xi^{m}, \xi^{-m(d^{n-1} + \ldots + d)}, \xi^{-m(d^{n-2} + \ldots + d)}, \ldots, \xi^{-m(d^{2} + d)}, \xi^{-md})
                                \end{equation}
                                be a singular point of $X_{0}$, where $\xi$ is a primitive root of unity of order $N.$ For every  $1 \le m \le N$  and every $\alpha\in\U$ let $p_{m,\al}$ be the corresponding singular point of $X_{\al}$ close to $p_{m}$. Consider the Jacobian matrix of $X_{\alpha}$ at the singular point $p_{m,\alpha} = (x_{1,\alpha}, x_{2,\alpha}, \ldots x_{n,\al}).$
The characteristic polynomial is given by
\[
P_{\alpha}(\lambda) = \det(\lambda I - DX_{\alpha}(p_{m,\alpha})) = \lambda^n + \sum_{i=1}^{n} \sigma_{i,m,\alpha} \lambda^{n-i}.
\]
Now, consider the map 
\begin{equation}   \label{e:Phi_star}    \begin{split}  \Phi^\circ_{m} : \U &\rightarrow \C^n,\\
                  \Phi^\circ_{m}(\alpha) &:= (\sigma_{1,m,\alpha}, \sigma_{2,m,\alpha}, \ldots ,\sigma_{n,m,\alpha}).
 \end{split}
\end{equation}
The main result of this  section   is
the following

\begin{theorem}\label{T:submersion}
 It  holds that $\rank(D\Phi^\circ_{m}(0)) = n$ for all $1 \leq m \leq N$, where
\[
D\Phi^\circ_{m}(0) = \begin{pmatrix}
\frac{\partial \sigma_{1,m,\alpha}}{\partial \alpha_1}(0)& \frac{\partial \sigma_{1,m,\alpha}}{\partial \alpha_2}(0) & \ldots & \frac{\partial \sigma_{1,m,\alpha}}{\partial \alpha_n}(0)\\
\frac{\partial \sigma_{2,m,\alpha}}{\partial \alpha_1}(0) & \frac{\partial \sigma_{2,m,\alpha}}{\partial \alpha_2}(0) & \ldots & \frac{\partial \sigma_{2,m,\alpha}}{\partial \alpha_n}(0)\\
\vdots & \vdots & \ddots & \vdots\\
\frac{\partial \sigma_{n,m,\alpha}}{\partial \alpha_1}(0) & \frac{\partial \sigma_{n,m,\alpha}}{\partial \alpha_2}(0) & \ldots & \frac{\partial \sigma_{n,m,\alpha}}{\partial \alpha_n}(0)\\
\end{pmatrix}.
\]
\end{theorem}

For the proof of the  proposition,  we need several lemmas.     The first lemma gives  an asymptotic  description of   the singular set $E_{\alpha}$ in terms of $\alpha.$

\begin{lemma}\label{L:asym-sing-sets} Let $x=(x_{1,\alpha}, x_{2,\alpha}, \ldots x_{n,\al})$ be a point in  $ E_{\alpha}.$ Then for all $1 \le i \le n,$ we have
\begin{equation}\label{e:x_i,alpha}
x_{i,\al} = x_{1,\al}^{\sum_{s=1}^{n+1-i} -d^{s}} + \alpha_{i} x_{1,\al}^{-d} + \sum_{j = i+1}^{n} \alpha_{j} d^{j-i} x_{1,\al}^{-d + \sum_{l = 0}^{j-i-1}-d^{n+1-i-l}} + O(|\alpha|^2).
\end{equation}
\end{lemma} 
\proof We will prove it by descending induction on $i\in[1,n].$

Since $x\in E_\alpha,$  we have  by \eqref{e:E_alpha} $x_{n,\alpha} x_{1,\alpha}^{d} = 1 + \alpha_{n}$, and hence
\[
x_{n,\alpha} = x_{1,\alpha}^{-d} + \alpha_{n}x_{1,\alpha}^{-d},
\]
which shows that   estimate \eqref{e:x_i,alpha} is true for $i=n.$ Now suppose that estimate \eqref{e:x_i,alpha} is true for $i=m$. Writing  
\[
x_{m,\al}^{d} = \left(x_{1,\al}^{\sum_{s=1}^{n+1-m} -d^{s}} + \alpha_{m} x_{1,\al}^{-d} + \sum_{j = m+1}^{n} \alpha_{j} d^{j-m} x_{1,\al}^{-d + \sum_{l = 0}^{j-m-1}-d^{n+1-m-l}} +      O(|\alpha|^2 )    \right)^{d},
\]
we get
\begin{multline*}
    x_{m,\al}^{d} = \left(x_{1,\al}^{\sum_{s=1}^{n+1-m} -d^{s}}\right)^{d} + \binom{d}{d-1} \left(x_{1,\al}^{\sum_{s=1}^{n+1-m} -d^{s}}\right)^{d-1} \\\left(\alpha_{m} x_{1,\al}^{-d} + \sum_{j = m+1}^{n} \alpha_{j} d^{j-m} x_{1,\al}^{-d + \sum_{l = 0}^{j-m-1}-d^{n+1-m-l}}\right) + O(|\alpha|^2).
\end{multline*}
Since $(d-1)(\sum_{s=1}^{n+1-m} -d^{s}) = \left(\sum_{s=2}^{n+1-(m-1)} -d^{s}\right) + \left(\sum_{s=1}^{n+1-m} d^{s}\right) = d - d^{n+1-(m-1)}$, we have
\begin{multline*}
    x_{m,\al}^{d} = x_{1,\al}^{\sum_{s=2}^{n+1-(m-1)} -d^{s}} + d x_{1,\al}^{d - d^{n+1-(m-1)}} \\\left(\alpha_{m} x_{1,\al}^{-d} + \sum_{j = m+1}^{n} \alpha_{j} d^{j-m} x_{1,\al}^{-d + \sum_{l = 0}^{j-m-1}-d^{n+1-m-l}}\right) +  O(|\alpha|^2 ) .
\end{multline*}
Therefore, we get
\begin{equation*}
    x_{m,\al}^{d} = x_{1,\al}^{\sum_{s=2}^{n+1-(m-1)} -d^{s}} + \sum_{j = (m-1)+1}^{n} \alpha_{j} d^{j-(m-1)} x_{1,\al}^{\sum_{l = 0}^{j-(m-1)-1}-d^{n+1-(m-1)-l}} +  O(|\alpha|^2 ) .
\end{equation*}
Now using the above estimate and the relation $x_{m-1,\al} x_{1,\al}^{d} = x_{m,\al}^{d} + \al_{m-1}$  obtained from \eqref{e:E_alpha}, it follows that
\begin{multline*}
x_{m-1,\al} = x_{1,\al}^{-d}(x_{m,\al}^{d} + \al_{m-1}) = x_{1,\al}^{\sum_{s=1}^{n+1-(m-1)} -d^{s}} + \al_{m-1} x_{1,\al}^{-d} \\+ \sum_{j = (m-1)+1}^{n} \alpha_{j} d^{j-(m-1)} x_{1,\al}^{-d + \sum_{l = 0}^{j-(m-1)-1}-d^{n+1-(m-1)-l}} +  O(|\alpha|^2 ) .
\end{multline*}
So estimate  \eqref{e:x_i,alpha} is true for $i=m-1,$ and the proof is thereby completed.  
\endproof

\medskip

\begin{lemma}\label{L:Char_poly}
 With the convention   $x_{0,\alpha} = 1$, we have that
\[
P_{\alpha}(\lambda) = (\lambda + x_{1,\alpha}^d)^{n} + \sum_{j=1}^{n}{\left[d^j (x_{1,\alpha} x_{2,\alpha} \ldots x_{j-1,\alpha})^{d-1} x_{j,\alpha}^{d} (\lambda + x_{1,\alpha}^d)^{n-j}\right]}.
\]
Moreover, it holds  for $1 \leq i \leq n$ that
\[
\sigma_{i,m,\al} =\binom{n}{n-i} x_{1,\al}^{id} + \sum_{j=1}^{i} \binom{n-j}{n-i} x_{1,\al}^{(i-j)d} (d^{j}(x_{1,\al} x_{2,\al} \ldots x_{j-1,\al})^{d-1} x_{j,\al}^{d}).
\]
\end{lemma}
\proof
Using the  explicit expression  \eqref{e:Jouanolou-new},  a straightforward computation gives 
  the following Jacobian matrix of $X_{\alpha}$ at the singular point $p_{m,\alpha} = (x_{1,\alpha}, x_{2,\alpha}, \ldots x_{n,\al})$
\[
DX_{\alpha}(p_{m,\alpha}) = \begin{pmatrix}
-(d+1)x_{1,\alpha}^{d} & dx_{2,\alpha}^{d-1} & 0 & 0 & \ldots & 0 & 0\\
-dx_{2,\alpha} x_{1,\alpha}^{d-1}& -x_{1,\alpha}^{d} & d x_{3,\alpha}^{d-1} & 0 & \ldots & 0 & 0\\
-dx_{3,\alpha} x_{1,\alpha}^{d-1}& 0 & -x_{1,\alpha}^d  & dx_{4,\alpha}^{d-1} & \ldots & 0 & 0\\
\vdots& \vdots & \vdots & \vdots & \ddots & \vdots & \vdots\\
-dx_{n-1,\alpha} x_{1,\alpha}^{d-1}& 0 & 0 & 0 & \ldots & -x_{1,\alpha}^{d} & d x_{n,\alpha}^{d-1}\\
-dx_{n,\alpha} x_{1,\alpha}^{d-1}& 0 & 0 & 0 & \ldots & 0 & -x_{1,\alpha}^{d}\\
\end{pmatrix}.
\]
Therefore, we get that 
\[
\lambda I - DX_{\alpha}(p_{m,\alpha}) = \begin{pmatrix}
\lambda +(d+1)x_{1,\alpha}^{d} & -dx_{2,\alpha}^{d-1} & 0 & 0 & \ldots & 0 & 0\\
dx_{2,\alpha} x_{1,\alpha}^{d-1}& \lambda + x_{1,\alpha}^{d} & -d x_{3,\alpha}^{d-1} & 0 & \ldots & 0 & 0\\
dx_{3,\alpha} x_{1,\alpha}^{d-1}& 0 & \lambda + x_{1,\alpha}^d & -dx_{4,\alpha}^{d-1} & \ldots & 0 & 0\\
\vdots& \vdots & \vdots & \vdots & \ddots & \vdots & \vdots\\
dx_{n-1,\alpha} x_{1,\alpha}^{d-1}& 0 & 0 & 0 & \ldots & \lambda +x_{1,\alpha}^{d} & -d x_{n,\alpha}^{d-1}\\
dx_{n,\alpha} x_{1,\alpha}^{d-1}& 0 & 0 & 0 & \ldots & 0 & \lambda +x_{1,\alpha}^{d} \\
\end{pmatrix}.
\]
 The Laplace expansion along the first  column  of the above  matrix gives  
 \begin{eqnarray*}
  P_\alpha(\lambda)&=&\det\big(\lambda I - DX_{\alpha}(p_{m,\alpha})\big)=  [\lambda +(d+1)x_{1,\alpha}^{d} ]( \lambda +x_{1,\alpha}^{d})^{n-1}\\
  &+&\sum_{j=2}^n (-1)^{j-1}d x_{j,\alpha}x_{1,\alpha}^{d-1}\det
  \begin{pmatrix}
 L_{j-1} &  \rvline &  0\\ 
\hline
  \ 0 & \rvline & U_{n-j} 
\end{pmatrix}
\end{eqnarray*}
where $L_{j-1}$  is  a lower triangular $(j-1)\times (j-1)$  matrix and $U_{n-j}$  is a upper triangular $(n-j)\times (n-j)$  matrix of the form 
\[ L_{j-1}=
 \begin{pmatrix}
 -dx_{2,\alpha}^{d-1} \\
 & \ddots &   \\
  && -dx_{j,\alpha}^{d-1}
\end{pmatrix}\quad\text{and}\quad 
 U_{n-j} =
 \begin{pmatrix}
 \lambda+x_{1,\alpha}^{d} \\
 & \ddots &   \\
  && \lambda+x_{1,\alpha}^{d}
\end{pmatrix}.
\]
Since we have
$$
\det
  \begin{pmatrix}
 L_{j-1} &  \rvline &  0\\ 
\hline
  \ 0 & \rvline & U_{n-j} 
\end{pmatrix}
=\det L_{j-1} \det U_{n-j}=(-1)^{j-1}d^{j-1}(x_{2,\alpha}\ldots x_{j,\alpha})^{d-1}(\lambda +x_{1,\alpha}^d)^{n-j},
$$
a calculation  gives the  first  assertion of the lemma.

The second assertion is an  immediate consequence of the first one. 
\endproof
\begin{lemma}\label{L:asym-sing-sets-II}  For each $2 \le j \le n,$  there exist constants $c_{1,j}, c_{2,j}, \ldots , c_{j-2, j}\in\C$ such that $|c_{1,j}|=|c_{2,j}|= \ldots =|c_{j-2, j}|=1$ and that for $x=(x_{1,\alpha}, x_{2,\alpha}, \ldots x_{n,\al})\in E_\alpha,$  we have 
\begin{equation}\label{e:prod_x_j,alpha}
(x_{1,\al} x_{2,\al} \ldots x_{j-1,\al})^{d-1} x_{j,\al}^{d} = x_{1,\al}^{jd} + \sum_{l=1}^{j-2} c_{l,j} \al_{l} - x_{1,\al}^{\left((j-1)d - \sum_{l=1}^{j-1} d^{(n+2-j+l)}\right)} \al_{j-1} +  O(|\alpha|^2 ) .
\end{equation}
\end{lemma}
\proof %This is trivially true for $j=1$. 
Consider for $j=2$,
\[
x_{1,\al}^{d-1} x_{2,\al}^{d} = x_{1,\al}^{d-1} (x_{1,\al}^{d+1} - \al_{1}) = x_{1,\al}^{2d} - x_{1,\al}^{d-1} \al_{1}.
\]
We use  \eqref{e:x_i,alpha} for $i=1$ and  take the power of order $d-1$ on both sides to get
\[
x_{1,\al}^{d-1} = x_{1,\al}^{d - d^{n+1}} +  O(|\alpha|) .
\]
Therefore,
\[
x_{1,\al}^{d-1} x_{2,\al}^{d} = x_{1,\al}^{2d} - x_{1,\al}^{d- d^{n+1}} \al_{1} +  O(|\alpha|^2 ) .
\]
Now suppose \eqref{e:prod_x_j,alpha} is true for $j=m$. Consider
\[
\begin{split}
    (x_{1,\al} x_{2,\al} \ldots x_{m,\al})^{d-1} x_{m+1,\al}^{d} &= (x_{1,\al} x_{2,\al} \ldots x_{m,\al})^{d-1} (x_{m,\al} x_{1,\al}^{d} - \al_{m})\\
    &= x_{1,\al}^{d}(x_{1,\al} x_{2,\al} \ldots x_{m-1,\al})^{d-1} x_{m,\al}^{d} - (x_{1,\al} x_{2,\al} \ldots x_{m,\al})^{d-1} \al_{m}.\\
\end{split}
\]
Inserting  \eqref{e:prod_x_j,alpha} for $j=m$ into the first term in the last line, we get
\begin{multline*}
    (x_{1,\al} x_{2,\al} \ldots x_{m,\al})^{d-1} x_{m+1,\al}^{d} = x_{1,\al}^{(m+1)d} + \sum_{l=1}^{m-2} x_{1,\al}^{d} c_{l,m} \al_{l} - x_{1,\al}^{\left((m-1)d -   \sum_{l=1}^{m-1} d^{(n+2-m+l)}  \right)} \al_{m-1} \\- (x_{1,\al} x_{2,\al} \ldots x_{m,\al})^{d-1} \al_{m} +  O(|\alpha|^2 ) .
\end{multline*}
Applying  Lemma \ref{L:asym-sing-sets} (estimate \eqref{e:x_i,alpha}) for $j = 1,2,\cdots, m$  yields that
\[
\begin{split}
(x_{1,\al} x_{2,\al} \ldots x_{m,\al})^{d-1} \al_{m} &= \left(x_{1,\al}^{-\left[\sum_{s=1}^{n} d^{s} + \sum_{s=1}^{n-1} d^{s} + \cdots + \sum_{s=1}^{n+1-m} d^{s}\right]} +  O(|\alpha| ) \right)^{(d-1)} \al_{m}\\
&= \left(x_{1,\al}^{-\big (m \sum_{s=1}^{n+1-m} d^{s} + \sum_{l=1}^{m-1} (m-l) d^{n+1-m+l} \big)} + 
O(|\alpha| ) \right)^{(d-1)} \al_{m} \\
&= \left(x_{1,\al}^{\big( m d  - \sum_{l=1}^m d^{n+1-m+l}  \big)} +  O(|\alpha|) \right) \al_{m} \\
&= x_{1,\al}^{\big (       m d  - \sum_{l=1}^m d^{n+1-m+l}                  \big)} \al_{m}+  O(|\alpha|^2 ) \\
\end{split}
\]
 Consider, for $1\le l \le m-2$,  the constants $c_{l,m+1} := x_{1,0}^{d} c_{l,m}$, and $c_{m-1,m+1} := -x_{1,0}^{\left(      m d  - \sum_{l=1}^m d^{n+1-m+l}     \right)}.$  Since by \eqref{e:p_m} we have $|x_{1,0}|=1,$
 it follows that $|c_{l,m+1}|= | c_{l,m}|$ for $1\le l \le m-2$, and $|c_{m-1,m+1}|=1.$ Moreover, we infer that
\[
(x_{1,\al} x_{2,\al} \ldots x_{m,\al})^{d-1} x_{m+1,\al}^{d} = x_{1,\al}^{(m+1)d} + \sum_{l=1}^{m-1} c_{l,m+1,\al} \al_{l} - x_{1,\al}^{\big (  m d  - \sum_{l=1}^m d^{n+1-m+l}  \big)} \al_{m}+  O(|\alpha|^2 ) .
\]
This proves \eqref{e:prod_x_j,alpha} for $j=m+1,$ and  hence the lemma follows.
\endproof

\medskip

Now we arrive  at the 
\proof[End of the proof of Theorem \ref{T:submersion}]

By  Lemma \ref{L:asym-sing-sets}  (estimate \eqref{e:x_i,alpha} for $i=1$), we have
\[
x_{1,\al} = x_{1,\al}^{\sum_{s=1}^{n} -d^{s}} + \alpha_{1} x_{1,\al}^{-d} + \sum_{j = 2}^{n} \alpha_{j} d^{j-1} x_{1,\al}^{-d + \sum_{l = 0}^{j-2}-d^{n-l}} +  O(|\alpha|^2 ) .
\]
Multiplying both sides by $x_{1,\al}^{-1}$,  we get
\begin{equation}\label{e:x_1,alpha}
    1 = x_{1,\al}^{-N} + \alpha_{1} x_{1,\al}^{-1-d} + \sum_{j = 2}^{n} \alpha_{j} d^{j-1} x_{1,\al}^{-1 -d + \sum_{l = 0}^{j-2}-d^{n-l}} +  O(|\alpha|^2 ) .
\end{equation}
 By Lemma \ref{L:asym-sing-sets-II}, we can find for $1 \le i \le n$, the constants $\tilde{c}_{l,i}\in\C$ which depend on the constants  $c_{l,i}$  such that
\begin{equation}\label{e:sigma}
\sigma_{i,m,\al} = \left(\sum_{j=0}^{i} \binom{n-j}{n-i} d^{j}\right)x_{1,\al}^{id} + \sum_{l=1}^{i-2} \tilde{c}_{l,i} \al_{l} - d^{i} x_{1,\al}^{\left((i-1)d - \sum_{l=1}^{i-1} d^{n+2-i+l}\right)} \al_{i-1} +  O(|\alpha|^2 ) .
\end{equation}
On the other hand,  we infer from \eqref{e:x_1,alpha} that
\[
x_{1,\al}^{-N} + \alpha_{1} x_{1,\al}^{-1-d} + \sum_{j = 2}^{n} \alpha_{j} d^{j-1} x_{1,\al}^{-1 -d + \sum_{l = 0}^{j-2}-d^{n-l}} +  O(|\alpha|^2 )  = 1.
\]
By differentiating the above expression implicitly with respect to  $\al_{j}$ for $1 \le j \le n$ and then
evaluating it  at $\al = 0$, we get
\begin{equation}\label{e:sing-derivative}
-N (x_{1,0}^{-N-1})\frac{d x_{1,\al}}{d \al_{j}}(0) + d^{j-1} x_{1,0}^{-1 -d + \sum_{l = 0}^{j-2}-d^{n-l}} = 0.
\end{equation}

For the sake of clarity we divide the remainder of the  proof into two cases according to the singularity $p_m,$
$1\leq m\leq N.$
%%%%%%%%%%%%%%%%%%%%%%%%%%%%%%%%%%%%%%%%%%%%%%%%%%%%%%
\subsection{The singularity $p_{N} = (1,1,\ldots, 1)$}
%%%%%%%%%%%%%%%%%%%%%%%%%%%%%%%%%%%%%%%%%%%%%%%%%%%%%%%
Now suppose that $p_{N} = (1,1, \ldots, 1) \in X_{0}$. Therefore $x_{1,0} = 1$, and hence  we infer from 
\eqref{e:sing-derivative} that 
\begin{equation*}
\frac{d x_{1,\al}}{d \al_{j}}(0) = \frac{d^{j-1}}{N}.
\end{equation*}
Set $A_{n,d,i} := \sum_{j=0}^{i} \binom{n-j}{n-i} d^{j}.$ The above  equality, combined with  \eqref{e:sigma}, implies that  \[
\frac{d \sigma_{i,N,\al}}{d \al_{j}}(0) = \begin{cases}
id (A_{n,d,i}) \frac{d^{j-1}}{N} + \tilde{c}_{j,i,0} 
& \text{for } 1 \le j < i-1 \\
id (A_{n,d,i}) \frac{d^{i-2}}{N} - d^{i}
& \text{for } j = i-1\\
id (A_{n,d,i}) \frac{d^{j-1}}{N}
& \text{for } i-1 < j \le n
\end{cases}.
\]
Let $v_{i}$ be the $i$-th row vector of the matrix $D\Phi^\circ_{D}(0)$, that is 
\[
\begin{split}
v_{i} &= \left(\frac{d \sigma_{i,\al}}{d \al_{1}}(0), \frac{d \sigma_{i,\al}}{d \al_{2}}(0), \ldots, \frac{d \sigma_{i,\al}}{d \al_{n}}(0)\right)\\
&= \frac{i A_{n,d,i}}{N} \left(d,d^2, \ldots, d^{n}\right) + (\tilde{c}_{1,i,0}, \tilde{c}_{2,i,0}, \ldots, \tilde{c}_{i-2,i,0}, -d^{i},0, \ldots, 0).
\end{split}
\]
Also, since $A_{n,d,1} = (n+d)$, we get $v_{1} = \frac{(n+d)}{N} \left(d,d^2, \ldots, d^{n}\right)$. Therefore for $2 \le i \le n$, we have
\[
v_{i} = \frac{i A_{n,d,i}}{(n+d)}v_{1} + (\tilde{c}_{1,i,0}, \tilde{c}_{2,i,0}, \ldots, \tilde{c}_{i-2,i,0}, -d^{i},0, \ldots, 0).
\]
Since the determinant is a multilinear map of row  vectors, we can replace $v_{i}$ for $2 \le i \le n$ by $\tilde{v}_{i}$ when calculating $\det(D\Phi^\circ_{N}(0))$, where
\[
\tilde{v}_{i} = (\tilde{c}_{1,i,0}, \tilde{c}_{2,i,0}, \ldots, \tilde{c}_{i-2,i,0}, -d^{i},0, \ldots, 0).
\]
Therefore, we get
\[
\det{(D\Phi^\circ_{N}(0))} = \frac{(n+d)}{N} \det{\begin{pmatrix}
d& d^2 & d^3 & d^4 & \ldots & d^{n-1} & d^n\\
-d^{2}& 0 & 0 & 0 & \ldots & 0 & 0\\
\tilde{c}_{1,3,0}& -d^{3} & 0 & 0 & \ldots & 0 & 0\\
\vdots& \vdots & \vdots & \vdots & \ddots & \vdots & \vdots\\
\tilde{c}_{1,n,0}& \tilde{c}_{2,n,0} & \tilde{c}_{3,n,0} & \tilde{c}_{4,n,0} & \ldots & -d^{n} & 0\\
\end{pmatrix}.}
\]
Since $\det{(v_1,v_2,v_3, \ldots, v_{n-1}, v_{n})} = (-1)^{n-1} \det{(v_{2}, v_{3}, v_{4}, \ldots, v_{n}, v_{1})}$, we get
\[
\begin{split}
\det{(D\Phi^\circ_{N}(0))} &= \frac{(-1)^{(n-1)}(n+d)}{N} \det{\begin{pmatrix}
-d^{2}& 0 & 0 & 0 & \ldots & 0 & 0\\
\tilde{c}_{1,3,0}& -d^{3} & 0 & 0 & \ldots & 0 & 0\\
\vdots& \vdots & \vdots & \vdots & \ddots & \vdots & \vdots\\
\tilde{c}_{1,n,0}& \tilde{c}_{2,n,0} & \tilde{c}_{3,n,0} & \tilde{c}_{4,n,0} & \ldots & -d^{n} & 0\\
d& d^2 & d^3 & d^4 & \ldots & d^{n-1} & d^n\\
\end{pmatrix}}\\
&= \frac{(-1)^{(n-1)}(n+d)}{N} (-1)^{(n-1)} d^{(\sum_{j=1}^{n} j)} d^{(n-1)}\\
&= \frac{(n+d)}{N} d^{\frac{n(n+1)}{2}} d^{n-1} = \frac{(n+d)}{N} d^{\frac{(n^2 + n + 2n -2)}{2}}.
\end{split}
\]

\medskip
\noindent Thus $\det{(D\Phi^\circ_{N}(0))} = \frac{(n+d)}{N} d^{\frac{(n^2 + 3n -2)}{2}}$, and therefore $\det{(D\Phi^\circ_{N}(0))} \neq 0$ for all $d \ge 1$ and for all $n \ge 2$.
Theorem \ref{T:submersion} is thereby completed  for the case of singularity $p_N.$
%%%%%%%%%%%%%%%%%%%%%%%%%%%%%%%%%%%%%%%%%%%%%%%%%%%%%%%%%%%%%%%%%%%%%%%%%%%%%%%
\subsection{Other singular points}
%%%%%%%%%%%%%%%%%%%%%%%%%%%%%%%%%%%%%%%%%%%%%%%%%%%%%%%%%%%%%%%%%%%%%%%%%%%%%%%%
Fix an arbitrary   singular point  $p_{m} $ of $X_{0}$ given in \eqref{e:p_m} for   $1 \leq m <N$. Write $\xi^{m} = \zeta$, and we have by \eqref{e:p_m} $$p_m = (\zeta, \zeta^{-(d^{n-1} + \ldots + d)}, \zeta^{-(d^{n-2} + \ldots + d)}, \ldots, \zeta^{-(d^{2} + d)}, \zeta^{-d}).$$ Therefore $x_{1,0} = \zeta.$ Since $\zeta^N=1$  we infer from \ref{e:sing-derivative}
\begin{equation*}
\frac{d x_{1,\al}}{d \al_{j}}(0) = \frac{d^{j-1}}{N}(\zeta)^{-d - \sum_{l=0}^{j-2} d^{n-l}}.
\end{equation*}
Set $A_{n,d,i} := \sum_{j=0}^{i} \binom{n-j}{n-i} d^{j}.$ The above  equality, combined with  \eqref{e:sigma}, implies that 
\[
\frac{d \sigma_{i,m,\al}}{d \al_{j}}(0) = \begin{cases}
\left(\frac{id (A_{n,d,i}) d^{j-1}}{N}\right)\zeta^{id-1} \zeta^{-d-\sum_{l=0}^{j-2} d^{n-l}} + \tilde{c}_{j,i,0} 
& \text{for } 1 \le j < i-1 \\
\left(\frac{id (A_{n,d,i}) d^{i-2}}{N}\right)\zeta^{id-1} \zeta^{-d-\sum_{l=0}^{i-3} d^{n-l}} - d^{i}\zeta^{c_{i}}
& \text{for } j = i-1\\
\left(\frac{id (A_{n,d,i}) d^{j-1}}{N}\right)\zeta^{id-1} \zeta^{-d-\sum_{l=0}^{j-2} d^{n-l}}
& \text{for } i-1 < j \le n
\end{cases},
\]
where the  $c_{i}$'s are integers for all $2 \le i \le n$.
Let $v_{i}$ be the $i$-th row vector of the matrix $D\Phi^\circ_{m}(0)$, that is 
\[
\begin{split}
v_{i} &= \left(\frac{d \sigma_{i,m,\al}}{d \al_{1}}(0), \frac{d \sigma_{i,m,\al}}{d \al_{2}}(0), \ldots, \frac{d \sigma_{i,m,\al}}{d \al_{n}}(0)\right)\\
&= \frac{(i A_{n,d,i} \zeta^{id-1})}{N} \left(d \zeta^{-d},d^2 \zeta^{-d-d^{n}}, d^{3} \zeta^{-d - d^{n-1} - d^{n}}, \ldots, d^{n} \zeta^{1-N}\right) \\ &+ (\tilde{c}_{1,i,0}, \tilde{c}_{2,i,0}, \ldots, \tilde{c}_{i-2,i,0}, -d^{i} \zeta^{c_{i}},0, \ldots, 0).
\end{split}
\]
Also, since $A_{n,d,1} = (n+d)$, we get $v_{1} = \frac{(n+d)\zeta^{(d-1)}}{N} \left(d \zeta^{-d},d^2 \zeta^{-d-d^{n}}, d^{3} \zeta^{-d - d^{n-1} - d^{n}}, \ldots, d^{n} \zeta^{1-N}\right)$. Therefore for $2 \le i \le n$, we have
\[
v_{i} = \frac{i A_{n,d,i}}{(n+d)} \zeta^{(i-1)d} v_{1} + (\tilde{c}_{1,i,0}, \tilde{c}_{2,i,0}, \ldots, \tilde{c}_{i-2,i,0}, -d^{i} \zeta^{c_{i}},0, \ldots, 0).
\]
 We  replace $v_{i}$ for $2 \le i \le n$ by $\tilde{v}_{i}$ when calculating $\det(D\Phi^\circ_{m}(0))$, where
\[
\tilde{v}_{i} = (\tilde{c}_{1,i,0}, \tilde{c}_{2,i,0}, \ldots, \tilde{c}_{i-2,i,0}, -d^{i} \zeta^{c_{i}},0, \ldots, 0).
\]
Therefore, $\det{(D\Phi^\circ_{m}(0))}$ is  equal to
\[
 \frac{(n+d) \zeta^{d-1}}{N} \det{\begin{pmatrix}
(d \zeta^{-d})& (d^2 \zeta^{-d-d^{n}}) & (d^3 \zeta^{-d-d^{n-1} - d^{n}}) & \ldots & (d^{n-1} \zeta^{1+d^{2}-N}) & (d^n \zeta^{1-N})\\
-d^{2} \zeta^{c_2}& 0 & 0 & \ldots & 0 & 0\\
\tilde{c}_{1,3,0}& -d^{3} \zeta^{c_3} & 0 & \ldots & 0 & 0\\
\vdots& \vdots & \vdots & \ddots & \vdots & \vdots\\
\tilde{c}_{1,n,0}& \tilde{c}_{2,n,0} & \tilde{c}_{3,n,0} & \ldots & -d^{n} \zeta^{c_{n}} & 0\\
\end{pmatrix}.}
\]
Since $\det{(v_1,v_2,v_3, \ldots, v_{n-1}, v_{n})} = (-1)^{n-1} \det{(v_{2}, v_{3}, v_{4}, \ldots, v_{n}, v_{1})}$, 
the above expression is  equal to
\[
\begin{split}
 &\quad \frac{(-1)^{(n-1)}(n+d) \zeta^{d-1}}{N} \det{\begin{pmatrix}
-d^{2} \zeta^{c_2}& 0 & 0 & \ldots & 0 & 0\\
\tilde{c}_{1,3,0}& -d^{3} \zeta^{c_3} & 0 & \ldots & 0 & 0\\
\vdots& \vdots & \vdots & \ddots & \vdots & \vdots\\
\tilde{c}_{1,n,0}& \tilde{c}_{2,n,0} & \tilde{c}_{3,n,0} & \ldots & -d^{n} \zeta^{c_{n}} & 0\\
(d \zeta^{-d})& (d^2 \zeta^{-d-d^{n}}) & (d^3 \zeta^{-d-d^{n-1} - d^{n}}) & \ldots & (d^{n-1} \zeta^{1+d^{2}-N}) & (d^n \zeta^{1-N})\\
\end{pmatrix}}\\
&= \frac{(-1)^{(n-1)}(n+d) \zeta^{d-1}}{N} (-1)^{(n-1)} d^{(\sum_{j=1}^{n} j)} d^{(n-1)} \zeta^{(1 + \sum_{i=2}^{n} c_{i})}\\
&= \frac{(n+d)}{N} d^{\frac{n(n+1)}{2}} d^{n-1} \zeta^{(d + \sum_{i=2}^{n} c_{i})} = \frac{(n+d)}{N} d^{\frac{(n^2 + n + 2n -2)}{2}} \zeta^{(d + \sum_{i=2}^{n} c_{i})}.
\end{split}
\]

\medskip
\noindent Thus $\det{(D\Phi^\circ_{m}(0))} = \frac{(n+d)}{N} d^{\frac{(n^2 + 3n -2)}{2}} \zeta^{(d + \sum_{i=2}^{n} c_{i})}$, and therefore $\det{(D\Phi^\circ_{m}(0))} \neq 0$ for all $d \ge 1$ and for all $n \ge 2$.
This completes the proof of  Theorem \ref{T:submersion}. 
\endproof

\begin{remark} \label{R:LinsNetoSoares}\rm
 It is tempting   to prove Theorem \ref{T:submersion} for the following $n$-dimensional family
 $$
Y_{\alpha} = \sum_{i=1}^{n-1}(x_{i+1}^d - x_{i}x_1^{d} + \alpha_{i}x_i)\frac{\partial}{\partial x_i} + (1 - x_n x_1^{d} + \alpha_{n}x_n)\frac{\partial}{\partial x_n} = X_{0} + \sum_{i=1}^{n} \alpha_{i}x_i \frac{\partial}{\partial x_{i}}.
 $$
  which is   very close to  the $1$-dimensional  family $X^d_\mu$  introduced by Lins Neto--Soares
 \cite[Theorem 1]{LinsNetoSoares}. However, the calculation of the corresponding characteristic polynomials of $Y_\alpha$  at its singular points   are  very complicated. Unfortunately, the authors are able to  prove  the analogue of  Theorem \ref{T:submersion} for the above family $(Y_\alpha)$ only for some small values of dimension $n.$
\end{remark}

%%%%%%%%%%%%%%%%%%%%%%%%%%%%%%%%%%%%%%%%%%%%%%%%
%%%%%%%%%%%%%%%%%%%%%%%%%%%%%%%%%%%%%%%%%%%%%%%%
 \section{Proof of the Main Theorem}\label{S:Main-Theorem }
%%%%%%%%%%%%%%%%%%%%%%%%%%%%%%%%%%%%%%%%%%%%%%%%
%%%%%%%%%%%%%%%%%%%%%%%%%%%%%%%%%%%%%%%%%%%%%%%%

%%%%%%%%%%%%%%%%%%%%%%%%%%%%%%%%%%%%%%%%%%%%%%%%%%%%%%%%%%%%%%%%%%%%%%%%%%%%%%%%%%%%%%%%%%%%%%%%%%
\subsection{No invariant algebraic curve for $\mathcal J_\alpha$ for small generic $\alpha\in\C^n$} \label{SS:No-invariant-curve}
%%%%%%%%%%%%%%%%%%%%%%%%%%%%%%%%%%%%%%%%%%%%%%%%%%%%%%%%%%%%%%%%%%%%%%%%%%%%%%%%%%%%%%%%%%%%%%%%%%

\begin{theorem} \label{T:no-inv-alg-curve} There is a  small  open neighborhood $\U$ of $0$ in $\C^n$ with  the  following properties:
\begin{enumerate}
 \item  If $n$ is  even, then 
for every $d\geq 2$ and $\alpha\in\U,$ the  foliation $\mathcal J_\alpha$ possesses no invariant algebraic curve.

\item  If $n$ is  odd, then for every  
for every $d\geq 2,$ there are  $$K(n,d):=d^{n-1}+d^{n-3}+\ldots+d^2+1$$
hyperplanes $\H^{n,d}_i$ in $\C^n$  such that   for every $\alpha\in\U\setminus  \bigcup_{i=1}^{K(n,d)} \H^{n,d}_i ,$ the  foliation $\mathcal J_\alpha$ possesses no invariant algebraic curve.
\end{enumerate}
\end{theorem}

Fix $n,d\in\N$ with $n\geq 2$ and $d\geq 1.$ 
Let $[a_1,\ldots,a_{n+1}]\in \PGL(n+1,\C)$ denote the class of the  matrix $\diag(a_1,\ldots,a_{n+1})$ and 
$\G(n,d)\subset \PGL(n+1,\C)$ be the subgroup consisting of elements   $[a_1,\ldots,a_{n+1}]\in \PGL(n+1,\C)$
which satisfy
\begin{equation*}
 {a_i\over a_{i+1}}=\big({a_{i+1}\over a_{i+2}}\big)^d,\qquad 1\leq i\leq n-1,\qquad {a_n\over a_{n+1}}=\big({a_{n+1}\over a_{1}}\big)^d.
\end{equation*}
It can be  checked that  the group $\G(n,d)$ is  cyclic of order $N=N(n,d)$ and is  generated by the class of
\begin{equation*}
\big [\xi, \xi^{-(d^{n-1}+\ldots+d)},\ldots, \xi^{-(d^2+d)},\xi^{-d},1\big ], 
\end{equation*}
where $\xi$ is a primitive  root of unity of order $N.$ Moreover, it can be shown that $\mathcal H_d$ acts freely and transitively on $\sing(\mathcal J_0).$

\begin{lemma} \label{L:pushforward} For every $\phi\in\G(n,d)$ and every $\alpha\in\C^n,$   we have   $\phi_*\mathcal J_\alpha=\mathcal J_{\tilde \alpha},$
where  $\tilde \alpha=\xi^{-d}\phi(\alpha).$
Here,   $\xi$  is a root of unity of order $N$ such that   for  all $x=(x_1,\ldots,x_n)\in\C^n,$
$$
\phi(x)=\big( \xi x_1, \xi^{-(d^{n-1}+\ldots+d)}x_2,\ldots, \xi^{-(d^2+d)}x_{n-1},\xi^{-d}x_n\big).$$
\end{lemma}
\proof
There are $\gamma_1,\ldots,\gamma_n\in\C^*$  such that
$$
\phi(x)=\big (\gamma^{-1}_1 x_1,\ldots, \gamma^{-1}_n x_n\big)\qquad\text{for}\qquad  x=(x_1,\ldots,x_n)\in\C^n.
$$
A straightforward computation  gives that
$$
\phi_* (X_0)=\sum\limits_{i=1}^{n-1}  (\gamma^d_{i+1} x_{i+1}^d -\gamma_i\gamma_1^d x_i x_1^d) \gamma_i^{-1}{\partial \over \partial x_i}
+ (1-\gamma_n\gamma_1^d  x_n x_1^d) \gamma_n^{-1}{\partial \over \partial x_n}.
$$
Since  
$$
\gamma_i^{-1}=\xi^{d+d^2+\ldots+d^{n-i+1}},\quad \gamma_{i+1}^d=\xi^{-(d^2+d^3+\ldots+d^{n-i+1})},\quad \gamma_1^d=\xi^d,\quad \gamma_n^{-1}=\xi^d,
$$
we  infer that
$$
\phi_*( X_\alpha)=\xi^d\big (  X_0+\sum\limits_{i=1}^n \tilde\alpha_i {\partial \over \partial x_i}\big),
$$
where 
$$
\tilde  \alpha=\big(\tilde\alpha_1,\ldots,\tilde\alpha_n\big)=\big( \xi^{d^2+\ldots+d^n} \alpha_1,
\xi^{d^2+\ldots+d^{n-1}} \alpha_2,\ldots, \xi^{d^2}\alpha_{n-1}, \alpha_n\big).
$$
Hence, the lemma  follows.
\endproof

\begin{lemma}\label{L:non-aligned-sing}
 There is a  small  open neighborhood $\U$ of $0$ in $\C^n$ with  the  following properties:
\begin{enumerate}
 \item  If $n$ is  even, then 
for every $d\geq 2$ and $\alpha\in\U,$  $\sing(\mathcal J_\alpha)$ does not have  $d+1$ points  aligned.

\item  If $n$ is  odd, then for every  
for every $d\geq 2,$ there are  $K(n,d)$  hyperplanes $\H^{n,d}_i$ in $\C^n$  such that   for every $\alpha\in\U\setminus  \bigcup_{i=1}^{K(n,d)} \H^{n,d}_i ,$  $\sing(\mathcal J_\alpha)$ does not have  $d+1$ points  aligned.
\end{enumerate}
\end{lemma}
\proof
Consider first the case where $n$ is  even. By \cite[Lemma 3.4, p. 665]{LinsNetoSoares}   
$\sing(\mathcal J_0)$ does not have  $d+1$ points  aligned.
By  Proposition \ref{P:Local-holo-charts}, we know that $\sing(\mathcal J_\alpha)$ depends holomorphically on $\alpha\in\U.$  By a geometric continuity argument, 
 we infer that when $\U$ is  small enough,  $\sing(\mathcal J_\alpha)$ does not have  $d+1$ points  aligned. This proves assertion (1).

 In the remainder of the  proof suppose  that $n$ is  odd.
Fix $\nu = (\nu_1, \nu_2, \ldots, \nu_n) \in \mathbb{C}^n\setminus \{0\}.$  For every $\mu\in\C$ consider the singular holomorphic foliation $\mathcal J_{\mu,\nu}$ on $\P^n$  which is generated on $\C^n$ by the  vector field
\begin{equation}\label{e:Jouanolou-new-mu}
X_{\mu} = \sum_{i=1}^{n-1}(x_{i+1}^d - x_{i}x_1^{d} + \mu\nu_{i})\frac{\partial}{\partial x_i} + (1 - x_n x_1^{d} + \mu\nu_{n})\frac{\partial}{\partial x_n} = X_{0} + \sum_{i=1}^{n} \mu\nu_{i} \frac{\partial}{\partial x_{i}}.
\end{equation}
By \eqref{e:Jouanolou-new}, we have  $\mathcal J_{\mu,\nu}=\mathcal J_{\mu\nu}.$

For  $\mu\in\C$ with $|\mu|\ll 1,$ let $p_{i,\mu}$ be the  singular point of $\mathcal J_{\mu,\nu}$ corresponding  to the singular point $p_{i,0}$ of $\mathcal J_{0,\nu}=\mathcal J_0$ for every $1\leq  i\leq N.$
Write
$$
p_{i,\mu}=\big(  x_{1,i,\mu},  x_{2,i,\mu},\ldots, x_{n,i,\mu}\big)\in \sing(\mathcal J_{\mu,\nu})\subset\C^n.
$$
We infer from  this and \eqref{e:Jouanolou-new-mu} that for every $1\leq  i\leq  N,$
\begin{equation}\label{e:x_i,mu}
 x^{-N}_{1,i,\mu}+\mu\big(\nu_1x^{-1-d}_{1,i,\mu}  + d\nu_2x^{-1-d-d^n}_{1,i,\mu}+\ldots+  d^{n-2}\nu_{n-1}x^{-1-d-\sum_{l=0}^{n-2} d^{n-l}}_{1,i,\mu}+ d^{n-1}\nu_nx^{-N}_{1,i,\mu} \big)+O(\mu^2)=1.
\end{equation}
By differentiating the above equation with respect to $\mu$ at $\mu=0,$ we get that
\begin{equation*}
 -Nx^{-N-1}_{1,i,0}\left.{d  x_{1,i,\mu}\over d\mu}\right|_{\mu=0}+\big(\sum\limits_{k=1}^n\nu_k d^{k-1} x^{-1-d-\sum_{l=0}^{k-2} d^{n-l}}_{1,i,0}\big)=0.
\end{equation*}
Hence, using  \eqref{e:N} we obtain that
\begin{equation*}
 \left.{d  x_{1,i,\mu}\over d\mu}\right|_{\mu=0}=\sum\limits_{k=1}^n\nu_k {d^{k-1}\over N} x^{N-d-\sum_{l=0}^{k-2} d^{n-l}}_{1,i,0}=\sum\limits_{k=1}^n\nu_k {d^{k-1}\over N} x^{-d+\sum_{l=0}^{n-k+1} d^{l}}_{1,i,0}.
\end{equation*}
Next, since 
$$
x_{n,i,\mu}=(1+\mu \nu_n)x_{1,i,\mu}^{-d},
$$
by differentiating the  equation with respect to $\mu,$ it follows that
\begin{eqnarray*}
 \left.{d  x_{n,i,\mu}\over d\mu}\right|_{\mu=0}&=&
 \nu_n x^{-d}_{1,i,0}+(-d)(x_{1,i,0})^{-d-1}\left.{d  x_{1,i,\mu}\over d\mu}\right|_{\mu=0}\\
 \\
 &=&\sum\limits_{k=1}^n-\nu_k {d^{k}\over N} x^{-1-2d+\sum_{l=0}^{n-k+1} d^{l}}_{1,i,0}+\nu_n x^{-d}_{1,i,0}.
\end{eqnarray*}

Let $\rho$ be a primitive root of unity of order $d+1,$ i.e. $\rho^{d+1}=1$ and $\rho^j\not=1$ for all $1\leq j\leq d.$ Consider the following $(d+1)$ points of $\sing(\mathcal J_0):$
 $$
 q_{j,0}:=\big( \rho^j,1,\rho^j,1,\ldots, 1, \rho^j \big),\qquad  0\leq j\leq d. 
 $$
 Then,  we know from  \cite[Lemma 3.4, p. 665]{LinsNetoSoares}  that the $d+1$ points $\phi (q_{0,0}),\ldots, \phi(q_{d,0})$ are aligned for every $\phi\in\G(n,d).$ In fact,  this describes the set of  all $(d+1)$-aligned points  of   $\sing(\mathcal J_0).$ The cardinality of this set is $ {N(n,d)\over d+1}=K(n,d)$ as   there   are 
 exactly $d+1$  elements $\phi\in \G(n,d)$ making  the set of  $d+1$ points 
 $\{q_{0,0},\ldots, q_{d,0}\}$ invariant.
 Therefore, we fix   a subset $\{ \phi_i:\ 1\leq i\leq K(n,d)\}\subset \G(n,d)$   such  that $\big \{  \{\phi_i (q_{0,0}),\ldots, \phi_i(q_{d,0})\}:\  1\leq i\leq K(n,d)\big \} $ enumerate all   sets of  all $(d+1)$-aligned points  of   $\sing(\mathcal J_0).$
 
For $0\leq j\leq d$ and $
 \mu\in\C$ with $|\mu|\ll 1,$ let $q_{j,\mu}$ be the corresponding   point of $q_{j,0}$  given by  Proposition \ref{P:Local-holo-charts}.   
 Suppose that $(q_{0,\mu},\ldots, q_{d,\mu})$ are aligned. So there is $t_2\in\C$ such that
$  t_2(q_{1,\mu}-q_{0,\mu})= q_{2,\mu}- q_{0,\mu}.   $ Write
$$
q_{j,\mu}=p_{i_j,\mu}=\big ( x_{1,i_j,\mu},\ldots, x_{n,i_j,\mu} \big).
$$
We deduce that
\begin{equation*}
 t_2(x_{1,i_1,\mu} -x_{1,i_0,\mu})=(x_{1,i_2,\mu} -x_{1,i_0,\mu})\quad\text{and}\quad 
 t_2(x_{n,i_1,\mu} -x_{n,i_0,\mu})=(x_{n,i_2,\mu} -x_{n,i_0,\mu}).
\end{equation*}
Cancelling out $t_2,$ this system of two equations implies that
\begin{equation}
 (x_{1,i_1,\mu} -x_{1,i_0,\mu})(x_{n,i_2,\mu} -x_{n,i_0,\mu})=(x_{1,i_2,\mu} -x_{1,i_0,\mu})(x_{n,i_1,\mu} -x_{n,i_0,\mu}).
\end{equation}
On the other hand, since $x_{1,i_j,0}=\rho^j$ and $\rho^{d+1}=1,$ we  have 
\begin{eqnarray*}
 \left.{d  x_{1,i_j,\mu}\over d\mu}\right|_{\mu=0}&=&\sum\limits_{k=1}^n\nu_k {d^{k-1}\over N} (\rho^j)^{-d+\sum_{l=0}^{n-k+1} d^l}=\sum\limits_{k=1}^n\nu_k {d^{k-1}\over N} (\rho^j)^{1+\sum_{l=0}^{n-k+1} d^l},\\
 \left.{d  x_{n,i_j,\mu}\over d\mu}\right|_{\mu=0}&=&\sum\limits_{k=1}^{n-1}-\nu_k {d^k\over N} (\rho^j)^{1+\sum_{l=0}^{n-k+1} d^l}-\nu_n {d^n\over N} (\rho^j)  +\nu_n(\rho^j).
\end{eqnarray*}
It follows that
\begin{equation}\label{e:condition-d}
 \sum_{k=1}^{n-1} \nu_kd^{k-1}P(k)=0,
\end{equation}
where 
\begin{equation}
 P(k):= \big (\rho^{1+\sum_{l=0}^{n-k+1} d^l}-1\big)(\rho^2-1)-(\rho-1)\big(\rho^{2(1+\sum_{l=0}^{n-k+1} d^l)}-1  \big).
\end{equation}
Since $n$ is  odd,  $n-k+1$ is  odd (resp. even) if and only if $k$ is  odd (resp. even).
As  we have
\begin{equation*}
 \sum_{l=0}^{n-k+1} d^l=\begin{cases}
                         (d+1)(d^{n-k}+d^{n-k-2}+\ldots+1),&\quad \text{if  $n-k+1$ is  odd;}\\
                          (d+1)(d^{n-k}+d^{n-k-2}+\ldots +d)+1,&\quad \text{if  $n-k+1$ is  even.}
                        \end{cases}
\end{equation*}
We see that
\begin{equation*}
 P(k)=\begin{cases}
       0,&\quad\text{if $k$ is  odd};\\
         (\rho^2-1)^2-(\rho-1)(\rho^4-1), &\quad\text{if $k$ is  even.}
      \end{cases}
\end{equation*}
So condition  \eqref{e:condition-d} becomes
\begin{equation}\label{e:condition-d-final}
 \left[ \sum_{k=1}^{n-1\over 2} \nu_{2k} d^{2k-1} \right] \left[(\rho^2-1)^2-(\rho-1)(\rho^4-1)\right]=0.
\end{equation}
The  expression in the second pair of brackets is  equal to $0$ if and only if $d=1.$ 
  Consider   the hyperplane  $\H^{n,d}$ given by 
$$
 \H^{n,d}:=\left\lbrace \nu=(\nu_1,\ldots,\nu_n)\in\C^n:\   \nu_2d+\nu_4 d^3+\ldots+\nu_{n-1}d^{n-2}=0 \right\rbrace.
$$
So  we have shown  that for every $\nu\not\in \H^{n,d},$  $\mathcal J_{\mu\nu}=\mathcal J_{\mu,\nu}$ does not have  $(q_{0,\mu},\ldots, q_{d,\mu})$ aligned  for $0<|\mu|\ll 1.$

By shrinking $\U$ if necessary, given $\alpha\in \U,$ 
let $q_{j,\alpha}$ be the corresponding   point of $q_{j,0}$  given by  Proposition \ref{P:Local-holo-charts} for every $0\leq j\leq d,$ and let $p_{i,\alpha}$ be the  singular point of $\mathcal J_{\alpha}$ corresponding  to the singular point $p_{i,0}$ of $\mathcal J_0$ for every $1\leq  i\leq N.$ Observe that  in \eqref{e:x_i,mu} the error  term   satisfies $|O(\mu^2)|\leq c|\mu|^2 $  for a constant
$c>0$  independent of 
$\nu\in\C^n$ with $\|\nu\|=1,$ using  the compactness of the unit sphere of $\C^n$ and  
a continuity argument.
Consequently, we may adapt  the arguments  from \eqref{e:x_i,mu} to \eqref{e:condition-d-final}
in order to   show  that if $\alpha\not\in\H^{n,d},$ then $(q_{0,\alpha},\ldots, q_{d,\alpha})$ are not aligned.  

%What  we have proved  so far is that  if  $(q_{0,t_k\alpha},\ldots, q_{d,t_k\alpha})$
%are aligned  for a  sequence   $(t_k)_{k=1}^\infty\subset \C\to 0,$ then 
%This  implies that $G(t_k\alpha)=0$ for all $k.$ Hence, $G(t\alpha)=0$ for all $t\in\C.$  

 Suppose in order to reach a contradiction that  $\sing(\mathcal J_\alpha)$ with $\alpha:=\mu\nu$ has  $d+1$ points  aligned.
By a geometric continuity argument as in the case of even $n,$ 
 we infer that when $\U$ is  small enough, 
 these  $d+1$ aligned points should be     the corresponding  points of   $\phi_i (q_{0,0}),\ldots, \phi_i(q_{d,0})$   given by  Proposition \ref{P:Local-holo-charts} for some $1\leq i\leq  K(n,d).$
 By Lemma \ref{L:pushforward},  $(\phi_i)_*\mathcal J_{\tilde \alpha}=\mathcal J_\alpha$
 with $\tilde \alpha:=\phi^{-1}_i(\xi^d_i\alpha),$ where $\xi_i$ is the root of unity of order $N$  defined by $\phi_i$ in Lemma \ref{L:pushforward}.  By the  above  paragraph, when  $\tilde\alpha\not\in \H^{n,d},$
 which is  equivalent to  the condition that $\xi^d_i\alpha\not\in \phi_i(\H^{n,d}),$ which is  in turn 
  equivalent to  $\alpha\not\in \phi_i(\H^{n,d}),$
 the above $d+1$  points are not aligned.  Therefore,    defining  for $1\leq i\leq  K(n,d),$ the hyperplane passing through  $0\in\C^n$:
 $$
 \H^{n,d}_i:= \phi_i(\H^{n,d}),
 $$
  we see that  assertion (2) of  the theorem holds.
\endproof

\proof[End of the proof of Theorem \ref{T:no-inv-alg-curve}]
We argue as  in the proof of \cite[Proposition 3.5]{LinsNetoSoares} using    Lemma \ref{L:non-aligned-sing}
instead of Lemma 3.4 therein.
\endproof
%%%%%%%%%%%%%%%%%%%%%%%%%%%%%%%%%%%%%%%%%%%%%%%%%%%%
\subsection{Proof of the Main Theorem}\label{SS:End}
%%%%%%%%%%%%%%%%%%%%%%%%%%%%%%%%%%%%%%%%%%%%%%%%%%%%%
If $d\geq 2,$ then using  Theorem \ref{T:no-inv-alg-curve}, we argue as in \cite[Section 4]{LinsNetoSoares} in order  to obtain  
a nonempty real Zariski  open set $\Uc^{\nd}_d(\P^n)\subset  \Fc_d(\P^n)$ such that for every $\Fc\in \Uc^{\nd}_d(\P^n),$
  all  the singularities of $\Fc$ are    non-degenerate and  $\Fc$ does not possess any invariant algebraic  curve
  (see also 
Remark \ref{R:LinsNetoSoares-bis}).  If $d=1,$ then  we set simply   $\Uc^{\nd}_d(\P^n):=\Fc^\nd_d(\P^n),$  see Definition \ref{D:F-nd}.

Fix an arbitrary   foliation   $\Fc_0\in \Fc^{\nd}_d(\P^n),$  and  set   $\sing(\Fc_0)=\{p_{1,0},\ldots, p_{N,0}\},$ where $p_{i,0}\not=p_{j,0}$ for $i\not=j.$ By Proposition  \ref{P:Local-holo-charts},
  there exist  connected open neighborhood $\mathcal U$ of $\Fc_0$ in $\Fc_d(\P^n),$   $U_1,\ldots, U_N$
 of $p_{1,0},\ldots, p_{N,0}$ respectively,  pairwise disjoint, and holomorphic maps $\varphi_j:\ \mathcal U\to  U_j,$
 for $1\leq j\leq N$  satisfying  the   properties (i)-(ii)-(iii) of the proposition.
Consider the maps
\begin{equation}\label{e:Phi_m}
\Phi_m:\  \mathcal U\to \C^n,\qquad m=1,\ldots, N,
\end{equation}
defined by
$$
\Phi_m(\Fc):=\big(\sigma_{1,m}(\Fc),\ldots, \sigma_{n,m}(\Fc)), 
$$
where $\sigma_{1,m}(\Fc),\ldots, \sigma_{n,m}(\Fc)$ are the coefficients of the  characteristic polynomial $P_{\Fc,m}(\lambda)$ of the Jacobian matrix of $X$ at the singular  point $\varphi_m(\Fc)\in U_m$ of $\Fc$, and $X$ is the polynomial  vector field in $\C^n$  representing $\Fc,$ that is,
$$
P_{\Fc,m}(\lambda) = \det(\lambda I - DX(\varphi_m(\Fc)) = \lambda^n + \sum_{i=1}^{n} \sigma_{i,m}(\Fc) \lambda^{n-i}.
$$
Consider the map 
\begin{eqnarray*} \theta:\C^n&\to&\C^n\\
 \lambda=(\lambda_1,\ldots,\lambda_n)&\mapsto &\big(\sigma_1(\lambda),\ldots,\sigma_n(\lambda)\big),
\end{eqnarray*}
where $\sigma_j,$ $j=1,\ldots,n,$ are the elementary symmetric functions of $\lambda_1,\ldots,\lambda_n.$
Consider  also  the set
$$
\mathcal D:= \left\lbrace   \lambda=(\lambda_1,\ldots,\lambda_n)\in\C^n:\ \exists\ 1\leq i <j\leq n:\ \lambda_i=\lambda_j \right\rbrace,
$$
which is  the union of ${n(n-1)\over 2}$ of complex hyperplanes in $\C^n,$ and the set $\Delta:=\theta(\mathcal D)$ which is a complex hypersurface   of $\C^n.$   Observe that
$$
\theta|_{\C^n\setminus \mathcal D}:\ \C^n\setminus \mathcal D\to \C^n\setminus \Delta
$$
is locally biholomorphic.

By shrinking  the  connected open neighborhood $\mathcal U$ of $\Fc_0$ in $\Fc_d(\P^n)$ and   $U_1,\ldots, U_N$  of $p_{1,0},\ldots, p_{N,0}$ if necessary, we may assume that there  exists a local  inverse $\theta_m$ of $\theta$ on $U_m$ for $1\leq m\leq N.$ Consider, for $1\leq m\leq N,$ the map
\begin{equation}\label{e:Psi_m}
 \Psi_m:\ \mathcal U\to \C^n\qquad \Psi_m:=\theta_m\circ\Phi_m.
\end{equation}
Then
\begin{equation}\label{e:Psi_m-bis}
 \Psi_m(\Fc)=(\lambda_1(\varphi_m(\Fc)),\ldots,  \lambda_n(\varphi_m(\Fc))  ),
\end{equation}
where $\lambda_j(\varphi_m(\Fc)),$ $j=1,\ldots, n,$ are the eigenvalues of the linear part of $\Fc$ at $\varphi_m(\Fc).$

We now arrive at the

\proof[End of the proof of  Theorem \ref{T:main}]
We will  prove  the following\\
\noindent  {\bf Fact.} {\it For every foliation $\Fc\in \Fc^\nd_d(\P^n),$ there exists a connected open neighborhood $\mathcal U=\mathcal U(\Fc)$ of $\Fc$ in  $\Fc^\nd_d(\P^n)$ and a  subset of $\Cc\subset \mathcal U$ of zero Lebesgue measure  such that  $  ( \mathcal U\setminus \Cc)\cap \Uc^{\nd}_d(\P^n)  \subset \Sc_d(\P^n).$  
}

Since we may find  a  countably many such  foliations  $\Fc$ such that $\mathcal U(\Fc)$ forms a countably  open  cover of $\Fc^\nd_d(\P^n),$
we  see that  the theorem  will follow from the above fact.

To prove  the above fact,  fix an arbitrary foliation $\Gc \in \Fc^\nd_d(\P^n).$  Since  $\Fc^\nd_d(\P^n)$
is  connected open,  we may find  a continuous  curve $\gamma:\ [0,1]\to  \Fc^\nd_d(\P^n)$ such that $\gamma(0)=\Fc_0:=\mathcal J^{n,d}$ the Jouanolou foliation  and  $\gamma(1)=\Gc.$

For a finite subdivision  $0=t_0 <t_1<\ldots <t_{k-1}<t_k=1$ of $[0,1],$ if we set $\Fc_\ell :=\gamma(t_\ell)$
then there is a connected open neighborhood $\mathcal U_\ell$ of $\Fc_\ell$ in $\Fc_d(\P^n),$ and for $1\leq m\leq N,$ the map $ \Phi_{m,\ell}$ given by \eqref{e:Phi_m}
and the map $\Psi_m$ given by \eqref{e:Psi_m}:
\begin{equation}\label{e:Psi_m-bisbis}
 \Phi_{m,\ell}:\  \mathcal U_\ell\to \C^n\quad\text{and}\quad \Psi_{m,\ell}:\ \mathcal U_\ell\to \C^n\qquad \Psi_{m,\ell}:=\theta_m\circ\Phi_{m,\ell}.
\end{equation}
 
By a compactness argument, we  may find such a subdivision  such that $\Fc_{\ell-1}\in  \mathcal U_\ell$ for $1\leq \ell\leq k$ and                                                 
  $\Fc_{\ell+1} \in  \mathcal U_\ell$ for $0\leq \ell< k.$                                               

 Fix a  $\delta>{n-2\over n}$ and set $\Ac:=(\C^n\setminus\Ac_\delta)\cup\Delta\subset\C^n.$ By  Proposition \ref{P:Arnold}, $\Ac$ is Borel set of zero $2n$-Lebesgue measure.
 
 It follows  from \eqref{e:Phi_star} and \eqref{e:Phi_m} that 
 \begin{equation}\Phi_{m,0}(\mathcal J_\alpha)=\Phi^\circ_m(\alpha)\qquad\text{for}\qquad \alpha\in\U. 
 \end{equation}
 This implies  that $\Im (D\Phi^\circ_m)(0)\subset \Im (D\Phi_{m,0})(\mathcal J_0)\subset \C^n.$
 On the other hand, by   Theorem \ref{T:submersion}, $ \Im (D\Phi^\circ_m)(0)= \C^n.$
 Therefore, we infer that   $\Im (D\Phi_{m,0})(\mathcal J_0)= \C^n,$ that is,  $\Phi_{m,0}$ is a submersion at the point $\mathcal J_0.$
 Hence, the critical set $$\mathcal C_{m,0}:=\{ \Fc\in \mathcal U_0:\ \Phi_{m,0} \quad \text{is not a submersion at}\quad \Fc\}$$
 is not equal to $\mathcal U_0.$  So  $\mathcal C_{m,0}$ is a proper analytic subset of $ \mathcal U_0$ since it is  defined  by the condition that  the determinant of every  minor $n\times n$ submatrix of the Jacobian matrix $D\Phi_{m,0}(\Fc)$ is equal to $0.$ By \eqref{e:Psi_m-bisbis}, $\mathcal C_{m,0}$ is also the critical set of 
 $ \Psi_{m,0}:\ \mathcal U_0\to \C^n .$   Consider the set
 \begin{equation*}
  \Cc_{0}:=\bigcup_{m=1}^N\Psi_{m,0}^{-1}(\Ac)\cup \mathcal C_{m,0}\subset \mathcal U_0.
 \end{equation*}
Since $\Psi_{m,0}$ is  a submersion outside the analytic set $\mathcal C_{m,0}$ and $\Ac$ is Borel set of zero $2n$-Lebesgue measure, it follows that $\Cc_{0}$  is  also a Borel set of zero $2M$-Lebesgue measure,
where  $M$ is  given by \eqref{e:M}.

Let $\Fc\in \mathcal U_0\setminus  \Cc_0.$  So $\Phi_{m,0}(\Fc)\not\in\Ac$ for all $1\leq m\leq N.$  By  Theorem \ref{T:Siegel}, 
  the vector  field  corresponding to $\Fc$  is  biholomorphically   equivalent  to its linear part in the neighborhood of  every  singular point $\varphi_m(\Fc).$  So $\Fc\in \Sc_d(\P^n).$  We have shown that
  $(\mathcal U_0\setminus  \Cc_0)  \cap \Uc^{\nd}_d(\P^n)\subset  \Sc_d(\P^n).$
  
  For $1\leq \ell\leq k$  consider the critical set $$\mathcal C_{m,\ell}:=\{ \Fc\in \mathcal U_\ell:\ \Phi_{m,\ell} \quad \text{is not a submersion at}\quad \Fc\}.$$
  Suppose that  we have  proved that $\mathcal C_{m,j}$ is an analytic set of $\mathcal U_j$ for every $0\leq j<\ell.$  We need to prove that  $\mathcal C_{m,\ell}$ is an analytic set of $\mathcal U_\ell.$
  It suffices  to show that $\mathcal C_{m,\ell}$ cannot be equal to the whole $\mathcal U_\ell.$
  Observe  that
  $$
  \mathcal C_{m,\ell-1} \cap \mathcal U_{\ell-1}\cap \mathcal U_\ell= \mathcal C_{m,\ell} \cap \mathcal U_{\ell-1}\cap \mathcal U_\ell.
  $$
  By the hypothesis of induction, the left-hand side is an analytic subset of  $\mathcal U_{\ell-1}\cap \mathcal U_\ell.$ This proves the last claim.
  So  $\mathcal C_{m,\ell}$ is a proper analytic subset of $ \mathcal U_\ell.$   Consider the set
 \begin{equation*}
  \Cc_{\ell}:=\bigcup_{m=1}^N\Psi_{m,\ell}^{-1}(\Ac)\cup \mathcal C_{m,\ell}\subset \mathcal U_\ell.
 \end{equation*}
Since $\Psi_{m,\ell}$ is  a submersion outside the analytic set $\mathcal C_{m,\ell}$ and $\Ac$ is Borel set of zero $2n$-Lebesgue measure, it follows that $\Cc_{\ell}$  is  also a Borel set of zero $2M$-Lebesgue measure. We can show as in the  case of $\mathcal U_0$ and $\Cc_0$ that
  $(\mathcal U_\ell\setminus  \Cc_\ell)\cap \Uc^{\nd}_d(\P^n)\subset  \Sc_d(\P^n).$
  
  When $\ell=k,$ setting $\mathcal U:=\mathcal U_k$ and $\Cc:=\Cc_k,$ we  obtain  that  $(\mathcal U\setminus  \Cc)\cap \Uc^{\nd}_d(\P^n)\subset  \Sc_d(\P^n),$ which completes the proof of the above fact.
\endproof

We conclude  the article with some remarks.
\begin{remark}\label{R:LinsNetoSoares-bis}
 \rm   Without   using  Theorem \ref{T:no-inv-alg-curve}
 we  can show  in Subsection \ref{SS:End} that there is     subset  $\Sc'_d(\P^n)$ of  $\Fc_d(\P^n)$  with full  Lebesgue  measure such that for every $\Fc\in\Sc'_d(\P^n), $ all singular points of  $\Fc$ are linearizable hyperbolic. This, combined  with Theorem \ref{T:generic-JNS} (see \cite[Theorem II]{LinsNetoSoares}),
 implies that
 Theorem \ref{T:main} holds for   $\Sc_d(\P^n):=\Sc'_d(\P^n)\cap  \Uc^\nd_d(\P^n).$
 
 Theorem \ref{T:no-inv-alg-curve} may be  regarded as the analogue of   $n$-dimensional parameter  of \cite[Theorem I]{LinsNetoSoares}.
 \end{remark}

 \begin{remark}\label{R:LorayRebelo} \rm By studying some groups of germs of holomorphic diffeomorphisms and the dynamics of such groups, Loray--Rebelo prove
 in \cite[Corollary B]{LorayRebelo}     the  following  elegant result:
 
 \noindent  {\bf Theorem.} {\it 
   There is     subset  $\Uc^{'}_d(\P^n)$ of  $\Fc_d(\P^n)$  with full  Lebesgue  measure such that  every $\Fc\in\Uc^{'}_d(\P^n) $  does not possess any invariant algebraic subvarieties of dimension $\geq 1.$
 }
 
 This result  was proved for $n=2$ in  Jouanolou \cite{Jouanolou} and for $n=3$ in Soares \cite{Soares}.
Combining this  result with   Theorem \ref{T:main}, we infer the  following  stronger version of  Theorem \ref{T:main} for
  $\Sc^{''}_d(\P^n):=  \Sc_d(\P^n)\cap \Uc^{'}_d(\P^n)$:
 \begin{corollary} There is     subset  $\Sc^{''}_d(\P^n)$ of  $\Fc_d(\P^n)$  with full  Lebesgue  measure such that for every $\Fc\in\Sc^{''}_d(\P^n), $ all singular points of  $\Fc$ are linearizable hyperbolic
 and that  $\Fc$ does not possess any invariant algebraic subvarieties of dimension $\geq 1.$
 \end{corollary}
\end{remark}

\begin{remark}
 \rm When $d=1,$  every foliation $\Fc\in\Fc_1(\P^n)$ possesses an invariant complex line. So the condition $d\geq 2$ in Theorem  \ref{T:main} in order to guarantee that there is no invariant algebraic  curve  is, in fact,  optimal.
\end{remark}

\begin{remark}\rm 
 It seems to be of interest  to  
 determine if  there is   an  open dense  set $\Uc_d(\P^n)$ of $\P^M$ (see \eqref{e:M}) such that  every $\Fc\in \Uc_d(\P^n)$ satisfies the conclusion of    Theorem \ref{T:main}.
 %We cannot hope that  the set $\Sc_d(\P^n)$ is  open dense  as in   Theorem \ref{T:generic-JNS}, since it is  not the case  %even in the local  context of germs of (local) holomorphic  vector fields  (see \cite{Arnold,IY}). 
\end{remark}

%%%%%%%%%%%%%%%%%%%%%%%%%%%%%%%%%%%%%%%%%%%%%%%%%

\end{document}